\documentstyle[11pt]{article}

\begin{document}
\begin{centering}
\huge\bf Braid group approach to the derivation of \\
universal R matrices
\vskip .6truecm
\normalsize\bf Feng Pan
\vskip .3truecm
\normalsize CCAST (World Lab.), P. O. Box 2375, Beijing 100080, P. R. China\\
\vskip .1truecm
Department of Physics, Liaoning Normal University, Dalian 116029, P. R. China
\vskip .4truecm
\bf Lianrong Dai
\vskip .3truecm
\normalsize Institute of High Energy Physics, Academia Sinica,
\vskip .1truecm 
P. O. Box 918-4-1, 100039 Beijing, P. R. China 
\vskip 1truecm
{\bf ABSTRACT}
\vskip .5truecm
\end{centering}
   A new method for deriving universal R matrices from braid
group representation is discussed. In this case, universal R
operators can be defined and expressed in terms of products of
braid group generators. The advantage of this method is that
matrix elements of R are rank independent, and leaves
multiplicity problem concerning coproducts of the corresponding
quantum groups untouched. As examples, R matrix elements of $[1]\times
[1]$, $[2]\times [2]$, $[1^{2}]\times [1^{2}]$, and $[21]\times [21]$
with multiplicity two for $A_{n}$, and $[1]\times [1]$ for $B_{n}$,
$C_{n}$, and $D_{n}$ type
quantum groups, which are related to Hecke algebra and Birman-Wenzl
algebra, respectively, are derived by using this method.
\newpage
\noindent{\large\bf 1. Introduction}
\vskip .3truecm
    Universal R matrices are solutions of spectral parameter-free
Yang-Baxter equations (YBE's). YBE's are of importance in both mathematics
and physics, such as statistical models$^{[1]}$, scattering matrices$^{[2]}$,
knot theory$^{[3]}$, conformal field theory (CFT)$^{[4]}$, and so on. Once the
parameter-free R matrices are known, the parameter dependent matrix \v{\bf R}(x)
can be obtained by using the so-called Baxterization procedure$^{[5-7]}$.
Up to now the derivations of the standard R matrices have been obtained
through the representation theory of quantum groups by many authors including
Drinfeld$^{[8]}$, Jimbo$^{[9-10]}$, Reshetikhin$^{[11]}$, and by taking
limit of statistical models$^{[12-14]}$ or by using the Witten's approach of
the link polynomials$^{[15-16]}$. There are also many other methods to
construct R matrices$^{[17-19]}$. Based on these methods, a various
class of R matrices have been obtained, which can easily be found
in the current mathematical physics literature. From these methods it can
easily be seen that R matrices are related either to the tensor products
of generators or to the CG coefficients of the corresponding quantum groups.
Thus knowledge of representation theory of quantum groups, such as coupling
coefficients, projection operators, and so on are very important in these
methods to construct the standard R matrices. In some cases, the
multiplicity problem will be involved in the coproducts of the corresponding
quantum groups, which is very complicated to solve. Secondly, R matrix
usually expressed in terms of CG matrix with summing over all the possible
resultant irreps of the corresponding quantum groups can only be derived for
specific $n$, e.g., of $A_{n}$, or $B_{n}$ at a time. When the $n$ increases,
the CG matrix will become very large. It will soon become intractable for
higher $n$ due to the drastic increase of the number of the CG coefficients
involved. It is also well known
that matrix representations of braid group generators can be constructed
by using the R matrix via
\vskip .3truecm
$$g_{i}=1\otimes 1\otimes \cdots \otimes 1\otimes R\otimes 1\otimes \cdots
\otimes 1, \eqno(1)$$
\vskip .3truecm
\noindent where R is in the $i$-th and $i+1$th spaces. It can  easily be
proved that the representations of the braid groups constructed in this way
are not irreducible in general. However, on the other hand, the R operator
can be regarded as deformed permutation operator which permutes two representations
of the corresponding quantum groups. From this point of view, R matrix is
representations of R operators in the uncoupled basis of the corresponding
quantum groups.
\vskip .3truecm
   In this paper, we will outline a new procedure for deriving standard solution
of the universal R matrices from representations of braid groups directly.
The R operators will be expressed in terms of products of the corresponding
braid group generators, which are acting on the vector product space of the
quantum groups. In this case, we only need CG coefficients of $[\lambda ]
\times [1]$ of the quantum groups. Such CG coefficients are very simple, and
of course, always multiplicity-free. Furthermore, the calculation is rank
independent. It opens up a new way to compute R matrix elements for
arbitrary $n$ once and for all, instead of one $n$ at a time.
As examples, R matrix elements of $[1]\times
[1]$, $[2]\times [2]$, $[1^{2}]\times [1^{2}]$, and $[21]\times [21]$
with multiplicity two for $A_{n}$, and $[1]\times [1]$ for $B_{n}$, $C_{n}$,
and $D_{n}$ type
quantum groups, which are related to Hecke algebra and Birman-Wenzl
algebra, respectively, are derived by using this method.
\vskip .5truecm
\noindent {\large\bf 2.~The R operators}
\vskip .3truecm
   Let $V^{[\lambda _{1}]}$ and $V^{[\lambda _{2}]}$ be spaces spanned by basis
vectors of irreducible representations of $[\lambda _{1}]$ and $[\lambda _{2}]$
of any quantum group. Then, the action of R is defined by
\vskip .3truecm
$$R(V^{[\lambda_{1}]}\otimes V^{[\lambda _{2}]})~\rightarrow
V^{[\lambda_{2}]}\otimes V^{[\lambda _{1}]}.\eqno(2)$$
\vskip .3truecm
\noindent We assume that maximum rank of $[\lambda_{1}]$ and $[\lambda_{2}]$ is
$f$. I. e., $[\lambda_{1}]$ and $[\lambda_{2}]$ cen be constructed by at most
$f${-}fold coproducts of rank-1 tensor operators of the corresponding
quantum group. For example, if $[\lambda_{1}]$ is of maximum rank, we can
write
\vskip .3truecm
$$T^{1}\otimes T^{2}\otimes \cdots \otimes T^{f}\rightarrow T^{[\lambda_{1}]},\eqno(3)$$
\vskip .3truecm
\noindent where $T^{i}$ with $i=1,~2,~\cdots ,f$ are vector operators in the
$i$-th space. In the A type quantum group case, in order to label the basis
of $T^{[\lambda_{1}]}$, one can assign a Weyl tableau $w^{[\lambda_{1}]}(\omega_{1}^{0})$
to $T^{[\lambda_{1}]}$, where $(\omega^{0}_{1})=(1,~2,~\cdots ,~f)$ is used to indicates
that the $f$  vectors are coupled to $[\lambda_{1}]$. Now we assume that the rank
of $[\lambda^{0}_{2}](\omega^{0}_{2},~2f-k+1,~2f-k+2,~\cdots ,~2f)$, where $(\omega^{0}_{2})=
(f+1,~f+2,~\cdots , 2f-k)$, while the indices $(2f-k+1,~2f-k+2,~\cdots ,~2f)$
are used to label the remaining scalars, i. e.,
\vskip .3truecm
$$T^{1}\otimes T^{2}\otimes \cdots \otimes T^{2f-k}\otimes 1^{2f-k+1}\otimes
\cdots\otimes 1^{2f}\rightarrow T^{[\lambda_{2}]}.\eqno(4)$$
\vskip .3truecm
\noindent Hence, the uncoupled basis vectors of $\vert w^{[\lambda_{1}]},
w^{[\lambda_{2}]}>$ can be written explicitly as
\vskip .3truecm
$$\vert w^{[\lambda_{1}]}(\omega_{1}^{0}),~
w^{[\lambda_{2}]}(\omega^{0}_{2},~2f-k+1,~2f-k+2,~\cdots ,~2f)>\eqno(5)$$
\vskip .3truecm
\noindent according to the braid group action, i. e., Eq. (5) can now be
understood as uncoupled basis vectors of braid group $B_{2f}$ under the R
operation with $R\{ (\omega_{1}^{0}),(\omega_{2}^{0},2f-k+1,2f-k+2,\cdots ,
2f)\}\rightarrow\{(\omega_{2}^{0}),2f-k+1, 2f-k+2,\cdots ,2f),(\omega_{1})\}$.
\vskip .3truecm
   Under this labelling scheme, we find the R operator can be expressed
in terms of braid group $B_{2f}$ generators.
\vskip .3truecm
$$R_{f=1}=g_{1},\eqno(6)$$
\vskip .3truecm
$$R_{f=2}=g_{2}g_{1}g_{3}g_{2}=g_{2}R_{1}g_{3}g_{2}.\eqno(7)$$
\vskip .3truecm
\noindent Through induction we finally obtain
\vskip .3truecm
$$R_{f}=g_{f}g_{f+1}g_{f+2}\cdots g_{2f-2}R_{f-1}g_{2f-1}g_{2f-2}\cdots
g_{f+2}g_{f+1}g_{f}.\eqno(8)$$
\vskip .3truecm
   (8) gives the universal R operator for fixed $f$
in braid group generator product form, where $f$ is the maximum rank of
irreps of the corresponding quantum groups.
The universality means that the operator given in (8) satisfies (2) for any
irrep of any type
of the corresponding quantum groups, which is just what Jimbo and Drinfeld
referred to$^{[8-10]}$.  The differences between this work and those of Jimbo
and Drinfeld are 1) the universal  R operator now is written in terms
of braid group generator product form; and 2) the R operator defined
in (8) is rank $f$ dependent. It seems that the R operator given in (8)
lose some universality. However, one can use it to compute all universal
R matrices because (8) is valid for arbitrary rank $f$ of irreps
in the tensor product space of the corresponding quantum groups given in (2).
While in Jimbo and Drinfeld's works
the R operator is written in terms of quantum double
basis of Hopf algebra, which is rank $f$ independent.
Thus (8) can be regarded as a braid group form
of the universal R operators formerly defined by Jimbo and Drinfeld.
\vskip .3cm
   Problem concerning the braid group realization for fixed type of the corresponding
quantum group has been studied by many works$^{[9-10,20]}$, from which one knows
that the braid group realization is Hecke algebra for the A type quantum groups,
is Birman-Wenzl algebra for the B, C, and D types, and is Kalfagianni algebra$^{[22]}$
for $G_{2}$. However, the problem still remains open for $F_{4}$, and the $E$ type
quantum groups. In the next section, we will outline a procedure for evaluating
the R matrices concerning A type quantum groups, and will also give a simple
example for the B, C, and D type cases.
\vskip .5cm
\noindent {\bf 3. Evaluation of R matrices}\\
\vskip .3truecm
   In this section, we will outline a procedure for evaluating R matrices.
We will consider Hecke algebra for the A type quantum groups and give an example
of Birman-Wenzl algebra for the B, C, and D type cases separately.
\vskip .3truecm
   The R operator is a braid group element, which can be expressed in terms of
braid group generators by using (8). R operates among the coordinate indices
$\{ 1,2,\cdots , 2f\}$. Any uncoupled basis vectors
$\vert w^{[\lambda_{1}]}(\omega_{1}^{0}),w^{[\lambda_{2}]} (\omega_{2}^{0},2f-k+1,2f-k+2,\cdots , 2f)>$
of any quantum group
can further be expanded in terms of uncoupled basis vectors of $2f-k${--}fold
basic representations, namely
\vskip .3truecm
$$\vert w^{[\lambda_{1}]}(\omega_{1}^{0}), w^{[\lambda_{2}]} (\omega_{2}^{0},2f-k+1,2f-k+2,\cdots , 2f)>=$$

$$\sum_{\omega}a_{\omega}Q_{\omega}\vert a_{1},a_{2},\cdots , a_{2f-k},1^{2k-k+1},1^{2f-k+2},\cdots , 1^{2f}>,\eqno(9)$$
\vskip .3truecm
\noindent where $a_{\omega}$ can be obtained by using the CG coefficients
for the coupling $((1\otimes )^{f})^{[\lambda_{1}]}((1\otimes )^{f-2k})^{[\lambda_{2}]}$
of the corresponding quantum group,
$\{ a_{1}, a_{2},\cdots , a_{2f-k}\}$ are the vector components of the
quantum group satisfying the normal ordering $a_{1}\leq a_{2}\leq \cdots
\leq a_{2f-k}$, and $Q_{\omega}$ is the left coset representatives in the
decomposition
\vskip .3truecm
$$B_{2f}=\sum_{\omega}\oplus Q_{\omega}(B_{1}\times B_{1}\times\cdots
\times B_{1}).\eqno(10)$$
\vskip .3truecm
\noindent For example, using the CG coefficients of $U_{q}(2)$ tabulated
in [21], we have
\vskip .3truecm
$$\vert aa, ab>=\sqrt{{q^{-1}\over{[2]}}}\vert a, a,a,b> +
\sqrt{{q\over{[2]}}}\vert a,a,b,a>$$

$$=(\sqrt{{q^{-1}\over{[2]}}}+\sqrt{{q\over{[2]}}}g_{3})\vert a,a,a,b>,\eqno(11a)$$
\vskip .3truecm
$$\vert ab, aa>=\sqrt{{q^{-1}\over{[2]}}}\vert a, b,a,a> +
\sqrt{{q\over{[2]}}}\vert b,a,a,a>$$

$$=(\sqrt{{q^{-1}\over{[2]}}}g_{2}g_{3}+\sqrt{{q\over{[2]}}}g_{1}g_{2}g_{3})\vert a,a,a,b>,\eqno(11b)$$
\vskip .3truecm
\noindent where $[x]$ is the q-number of x. The vector space indices are
arranged in natural order, e.g.,
\vskip .3truecm
$$\vert a,b,c>\sim T^{1}_{a}T^{2}_{b}T^{3}_{c},\eqno(12)$$
\vskip .3truecm
\noindent and the uncoupled basis vectors $Q_{\omega}\vert a_{1},a_{2},\cdots ,
a_{2f-k},1^{2f-k+1},1^{2f-k+2},\cdots ,1^{2f}>$ with different $\omega$ and
$a_{i}$'s are orthonormal.
\vskip .3truecm
$$< a_{1}^{\prime}, a_{2}^{\prime},\cdots a_{2f-k}^{\prime},
1^{2f-k+1},1^{2f-k+2},\cdots , 1^{2f}\vert Q^{\dagger}_{\omega^{\prime}}Q_{\omega}\vert
a_{1}, a_{2},\cdots a_{2f-k},1^{2f-k+1},1^{2f-k+2},\cdots , 1^{2f}>=$$

$$\delta_{\omega^{\prime}\omega}\prod\delta_{a_{i}^{\prime}a_{i}}.\eqno(13)$$
\vskip .3truecm
\noindent I. e., we use the orthogonal uncoupled basis of the
quantum group. In this case, it can be proved that the braid group parameters,
e. g., $q$, should
be real, otherwise (13) will no longer be valid. Equivalently, we have
used the following star operation
\vskip .3truecm
$$g_{i}^{\dagger}=g_{i},~~~{\rm for}~~~i=1,2,\cdots ,2f-1.\eqno(14)$$
\vskip .3truecm
\noindent However, results for generic braid group parameters can be
obtained through analytical continuation, i. e., the final
results are valid for generic parameters as well.
\vskip .3truecm
   The action of $g_{i}$ on the basis vectors $\vert a_{1},a_{2},\cdots ,
a_{2f-k}, 1^{2f-k+1},1^{2f-k+2},\cdots ,1^{2f}>$ is given by the following
rules
\vskip .3truecm
$$g_{i}\vert a_{1},a_{2},\cdots a_{2f-k},1^{2f-k+1},1^{2f-k+2},\cdots ,1^{2f}>$$

$$=q\vert a_{1},a_{2},\cdots a_{2f-k},1^{2f-k+1},1^{2f-k+2},\cdots ,1^{2f}>\eqno(15a)$$
\vskip .3truecm
\noindent if the components $a_{i}$ and $a_{i+1}$ are the same. This rule
can be proved by using the symmetrization method outlined in [21]. While
\vskip .3truecm
$$g_{i}\vert a_{1},a_{2},\cdots a_{2f-k},1^{2f-k+1},1^{2f-k+2},\cdots ,1^{2f}>$$

$$=\vert a_{1},\cdots , a_{i+1},a_{i},\cdots a_{2f-k},1^{2f-k+1},1^{2f-k+2},\cdots ,1^{2f}>\eqno(15b)$$
\vskip .3truecm
\noindent if the components $a_{i}$ and $a_{i+1}$ are different. It should be
noted that because of the property of braid groups, we should always write
the uncoupled basis vectors in the operator form. For example, in  Hecke
algebra case
\vskip .3truecm
$$g_{1}\vert a,b>=\vert b,a>,\eqno(16a)$$

$$g_{1}\vert b,a>=g^{2}_{i}\vert a,b>=(q-q^{-1})\vert b,a>+\vert a,b>,\eqno(16b)$$
\vskip .3truecm
\noindent otherwise the notation $\vert b,a>$ is rather confusing in the practical
computation.
\vskip .3truecm
   Hence, the action of the operator R on the uncoupled basis vectors of
quantum groups is well defined. Using Eqs. (15) and defining relations among
braid group generators, we can derive the R matrix elements. In the following,
we will outline how to use this procedure to derive R matrix elements. Firstly,
we will discuss Hecke algebra case.  Then, we will show you a simple example
in the Birman-Wenzl algebra case.
\vskip .5truecm
\noindent {\bf (1) Hecke algebra case}\\
\vskip .3truecm
   The Hecke algebra $H_{f}(q)$ is generated by $f-1$ elements
$g_{1},g_{2},\cdots , g_{f-1}$, which satisfy the following well-known
braid relations
\vskip .3truecm
$$g_{i}g_{i+1}g_{i}=g_{i+1}g_{i}g_{i+1},\eqno(17a)$$

$$g_{i}g_{j}=g_{j}g_{i},~~{\rm for \ }~~\vert i-j\vert\geq 2,\eqno(17b)$$

$$(g_{i})^{2}=g_{i}(q-q^{-1})+1.\eqno(17c)$$
\vskip .3truecm
   The braid group elements R, which can be expressed in terms of
braid group generators by using (8), are operating among the vector
space indices $\{ 1,2,\cdots ,2f\}$. Any uncoupled irreducible basis vectors
$\vert w^{[\lambda_{1}]}(\omega_{1}^{0}),w^{[\lambda_{2}]}(\omega_{2}^{0},
2f-k+1,2f-k+2,\cdots ,2f)>$ of $U_{q}(n)$ can further be expanded in
terms of uncoupled basis vectors of $2f-k${--}fold basic
representations as given by (9). In the following, we restrict ourselves
to $[\lambda_{1}]=[\lambda_{2}]$, and use $[2]\times [2]$ as an example.
\vskip .3truecm
\noindent {\bf Step 1.} Write out all the uncoupled basis vectors of the
corresponding quantum group.
\vskip .3truecm
In the $[2]\times [2]$ $U_{q}(n)$ case, these are
\vskip .3truecm
$$\vert ii,ii>=\vert i,i,i,i>,\eqno(18a)$$

$$\vert ij, ij>=(\sqrt{{q^{-1}\over{[2]}}}\vert i,j,1,1>+
\sqrt{{q\over{[2]}}}\vert j,i,1,1>)(\sqrt{{q^{-1}\over{[2]}}}\vert 1,1,i,j>+
\sqrt{{q\over{[2]}}}\vert 1,1,j,i>)$$

$$=A_{12}A_{34}g_{2}\vert i,i,j,j>,~~{\rm for\ }i< j, \eqno(18b)$$
\vskip .3truecm
\noindent where
\vskip .3truecm
$$A_{12}=\sqrt{1\over{[2]}}(q^{-1/2}+q^{1/2}g_{1}),\eqno(19a)$$

$$A_{34}=\sqrt{1\over{[2]}}(q^{-1/2}+q^{1/2}g_{3}),\eqno(19b)$$
\vskip .3truecm
\noindent where the CG coefficients of $[1]\times[1]\downarrow [2]$ of
$U_{q}(n)$ have been used. Similarly, we have
\vskip .3truecm
$$\vert ij,kl>=A_{12}A_{34}\vert i,j,k,l>, {\rm for\ }~~ i< j< k< l,\eqno(19c)$$
\vskip .3truecm
\noindent while
\vskip .3truecm
$$\vert kl,ij>=R\vert ij,kl>=g_{2}g_{3}g_{1}g_{2}\vert ij,kl>.\eqno(19d)$$
\vskip .3truecm
\noindent  All the other basis vectors can thus be written out similarly.
\vskip .3truecm
\noindent {\bf Step 2.} Derive algebraic relations among R, $A_{12}$, $A_{34}$.
\vskip .3truecm
  It can be proved that
\vskip .3truecm
$${R}A_{12}=A_{34}R,\eqno(20a)$$

$${R}A_{34}=A_{12}R,\eqno(20b).$$
\vskip .3truecm
\noindent Thus, we obtain
\vskip .3truecm
$${R}A_{12}A_{34}=A_{12}A_{34}R.\eqno(21)$$
\vskip .3truecm
\noindent (21) is very useful in the practical computation. We also need
the following quadratic equation of R
\vskip .3truecm
$$R^{2}=(q-q^{-1})g_{1}g_{3}R+(q-q^{-1})^{2}R+(q-q^{-1})g_{3}
g_{1}g_{2}g_{1}g_{3}$$

$$+(q-q^{-1})g_{2}g_{1}g_{2}+(q-q^{-1})g_{2}g_{3}g_{2}+(q-q^{-1})g_{2}+1.\eqno(22)$$
\vskip .3truecm
\noindent {\bf Step 3.} Applying R on all uncoupled basis vectors obtained
in step 1 and using the algebraic relations  derived in step 2 and
Eq. (15), we thus obtain
all the R matrix elements in this step. For example,
in $[2]\times [2]$ of $U_{q}(n)$ case, we have
\vskip .3truecm
$$R\vert ii,ii>=g_{2}g_{1}g_{3}g_{2}\vert i,i,i,i>=q^{4}\vert ii,ii>,\eqno(23)$$

$$R\vert ij,ij> ={R}A_{12}A_{34}\vert i,j,i,j>={R}A_{12}A_{34}g_{2}\vert i,i,j,j>$$

$$=A_{12}A_{34}g_{2}g_{3}g_{1}g_{2}^{2}\vert i,i,j,j>= $$

$$A_{12}A_{34}g_{2}g_{3}g_{1}g_{2}(q-q^{-1})\vert i,i,j,j>+
A_{12}A_{34}g_{2}g_{3}g_{1}\vert i,i,j,j>$$

$$=A_{12}A_{34}(q-q^{-1})\vert j,j,i,i>+q^{2}A_{12}A_{34}\vert i,j,i,j>.\eqno(24)$$
\vskip .3truecm
\noindent Using the following relation
\vskip .3truecm
$$A_{12}\vert i,i> = q^{1/2}[2]^{1/2}\vert i,i>,\eqno(25)$$
\vskip .3truecm
\noindent we obtain
\vskip .3truecm
$$R\vert ij,ij>=(q^{3}-q^{-1})\vert jj,ii>+q^{2}\vert ij,ij>,~~{\rm for\ } i< j\eqno(26)$$
\vskip .3truecm
\noindent Similarly, we have
\vskip .3truecm
$$R\vert ij,kl> =\vert kl,ij>, ~~{\rm for\ } i< j< k< l,\eqno(27a)$$
\vskip .3truecm
$$R\vert kl,ij> =R^{2}\vert ij,kl>=$$

$$(q-q^{-1})^{2}(q^{2}+1)\vert kl,ij>+(q^{3}-q^{-1})\vert jl,ik>+(q^{2}-1)\vert jk,il>$$

$$+(q^{2}-1)\vert il,jk>+(q-q^{-1})\vert ik,jl>+\vert ij,lk>,
{\rm for \ }~~ i< j< k< l.\eqno(27b)$$
\vskip .3truecm
\noindent All the other R matrix elements can thus be derived by using
this method. In the next section, we list all the R matrix elements for
$[1]\times [1]$, $[2]\times [2]$, $[1^{2}]\times [1^{2}]$, and some of $[21]
\times [21]$ for $U_{q}(n)$.
\vskip .3truecm
\noindent {\bf (2) A simple example of Birman-Wenzl algebra case}\\
\vskip .3truecm
   Birman-Wenzl algebra $C_{f}(r,q)$ is generated by $\{ g_{i}; i=1,2,\cdots ,
f-1\}$, which satisfy the braid group relations (17a), (17b) with constraint
\vskip .3truecm
$$(g_{i}-r^{-1})(g_{i}-q)(g_{i}+q^{-1})=0.\eqno(28)$$
\vskip .3truecm
   One can also introduce $f-1$ auxiliary generators $\{ e_{i}; i=1,2,\cdots ,
f-1\}$ with
\vskip .3truecm
$$e_{i}=1-{g_{i}-g_{i}^{-1}\over{q-q^{-1}}}.\eqno(29)$$
\vskip .3truecm
\noindent The following relations are helpful in practical computation.
\vskip .3truecm
$$e_{i}g_{i}=r^{-1}e_{i},\eqno(30a)$$
\vskip .3truecm
$$e^{2}_{i}=xe_{i},\eqno(30b)$$
\noindent
   Using the above defined relations, one obtains
\vskip .3truecm
$$g_{i}^{2}=(q-q^{-1})(g_{i}-r^{-1}e_{i})+1.\eqno(31)$$
\vskip .3truecm
   We consider $B_{n}$, $C_{n}$, and $D_{n}$ cases with irreps $[1]\times [1]$.
In this case $r=q^{2n}$ for $B_{n}$, $r=q^{2n-1}$ for $D_{n}$, and $r=q^{-2n-1}$ for $C_{n}$.
The procedure for evaluating the R matrix elements in this case is similar
to that for those of Hecke algebra case. Firstly, we write out all the uncoupled
basis vectors of $[1]\times [1]$. These are $\vert \mu , \nu>$ for $-n\leq \mu ,
\nu \leq n$. The R operator is  $g_{1}$. However, the $e_{1}$ operator is
so-called q-deformed trace contraction operator defined by
\vskip .3truecm
$$e_{1}\vert\mu ,-\nu >=(-)^{\mu}\delta_{\mu\nu}\sum_{k}(-)^{k}q^{k}\vert
k,~-k>\eqno(32)$$
\vskip .3truecm
\noindent when it is applied to the uncoupled basis vectors of the corresponding
quantum groups. Using (32), and assuming that the uncoupled basis vectors
are orthogonal, one can prove that
\vskip .3truecm
$$g_{1}\vert\mu ,-\mu>=r^{-1}\vert  -\mu , \mu>~~{\rm for\  }~~ \mu\neq 0.\eqno(33a)$$
\vskip .3truecm
\noindent If both $\mu$ and $\nu$ are nonzero, we can use (15) to derive
the R matrix elements. These are
\vskip .3truecm
$$g_{1}\vert\mu ,\mu >=q\vert\mu , \mu>,\eqno(33b)$$
\vskip .3truecm
$$g_{1}\vert\mu ,\nu >=\vert\nu ,\mu>,~~{\rm for\ }~~\mu>\nu ,~\mu\neq -\nu ,\eqno(33c)$$
\vskip .3truecm
\noindent Then, using (31), we get
\vskip .3truecm
$$g_{1}\vert\nu ,\mu >=(q-q^{-1})\vert\nu ,\mu>+\vert\mu , \nu>,~~{\rm for\ }~~\mu>\nu ,~
\mu\neq -\nu ,\eqno(33d)$$
\vskip .3truecm
$$g_{1}\vert -\mu ,\mu >=(q^{\mu -1}-q^{\mu +1}+r)\vert\mu ,-\mu >+
(q-q^{-1}-q^{1-\mu}+q^{-\mu-1})\vert -\mu ,\mu>$$

$$-(q-q^{-1})\sum_{k\neq\mu ,-\mu}q^{k}(-)^{\mu +k}\vert k,-k>,\eqno(33e)$$
\vskip .3truecm
In the $B_{n}$ case, we need uncoupled basis vectors $\vert 0,0>$.
Using (30), and (32), we obtain
\vskip .3truecm
$$g_{1}\vert 0,0>=\sum_{\mu >0}(-)^{\mu}a_{\mu}\vert -\mu ,\mu>+
\sum_{\mu >0}(-)^{\mu}b_{\mu}\vert\mu , -\mu > +c_{0}\vert 0, 0>,\eqno(33f)$$
\vskip .3truecm
\noindent where
$$a_{\mu}=q^{-\mu}r^{-1}-q^{-\mu +1}+q^{-\mu -1}+q^{1-2\mu}-q^{-2\mu -1}
-q^{\mu}r^{-1}$$

$$+{q-q^{-1}\over{q-1}}(q^{-\mu}-q^{-n-\mu}-q^{-2\mu +1}+q^{-2\mu }),\eqno(34a)$$
\vskip .3truecm
$$b_{\mu}=r^{-1}q^{\mu}-q^{-1}+q-rq^{-\mu}+{q-q^{-1}\over{q-1}}(q^{\mu}-q^{-n+\mu}
-q+1),\eqno(34b)$$
\vskip .3truecm
$$c_{0}=r^{-1}+{q-q^{-n+1}-q^{-1}+q^{-n-1}\over{q-1}}.\eqno(34c)$$
\vskip .3truecm
\noindent Using (31), we can prove that the basis vectors $\vert -\mu ,\mu>$, and
$\vert 0,0>$ are not orthonormal. The normalized basis vectors of $[1]\times [1]$
of B, C, and D cases are
\vskip .3truecm
$$\vert \mu ,\mu ) =\vert\mu ,\mu>,$$
\vskip .3truecm
$$\vert \mu ,\nu ) =\vert\mu ,\nu>~~{\rm for\ } \mu>\nu,\mu\neq -\nu,$$
\vskip .3truecm
$$\vert \nu ,\mu ) =\vert\nu ,\mu>~~{\rm for\ } \mu>\nu, $$
\vskip .3truecm
$$\vert\mu ,-\mu) =\vert \mu ,-\mu>,~~{\rm for \ }~~\mu >0$$
\vskip .3truecm
$$\vert -\mu ,\mu )={1\over{N_{\mu}}}\vert -\mu ,\mu>,~~{\rm for \ }\mu >0$$
\vskip  .3truecm
$$\vert 0,0 ) ={1\over{N_{0}}}\vert 0 ,0>,\eqno(35a) $$
\vskip .3truecm
\noindent where
\vskip .3truecm
$$N_{\mu} = \sqrt{r^{2}-rq^{\mu}(q-q^{-1})},$$

$$N_{0}^{2} ={\sum_{\mu >0}\{ a_{\mu}^{2}r(r-q^{\mu +1}+q^{\mu -1})+b^{2}_{\mu}\}
\over{1+(q-q^{-1})(c_{0}-r^{-1})+c^{2}_{0}}}.\eqno(35b)$$
\vskip .3truecm
   If one knows the CG coefficients for $[\lambda ]\times [1]$ for the corresponding
quantum groups, one can also derive other R matrix elements as have been done so
in the $A_{n}$ case.
\vskip.5truecm
\noindent {\bf 4. Some R matrix elements of $U_{q}(n)$ case}
\vskip .3truecm
\noindent (1) $[1]\times [1]$ irrep.
\vskip .3truecm In this case $R=g_{1}$. We have
\vskip .3truecm
$$R\vert i,i>=q\vert i,i>,$$

$$R\vert i,j> =\vert j,i>,$$

$$R\vert j,i>=(q-q^{-1})\vert j,i>+\vert i,j>. \eqno(36)$$
\vskip .3truecm
\noindent (2) $[2]\times [2]$ irrep.
\vskip .3truecm In this case $R=g_{2}g_{1}g_{3}g_{2}$. We have
\vskip .3truecm
$$R\vert ii,ii>=q^{4}\vert ii,ii>,$$

$$R\vert ii,ij>=q^{2}\vert ij,ii>,$$

$$R\vert ij,jj>=q^{2}\vert jj,ij>,$$

$$R\vert ii,jj>=\vert jj,ii>,$$

$$R\vert ij,kl>=\vert kl,ij>,$$

$$R\vert ij,kk>=\vert kk,ij>,$$

$$R\vert ik,jj>=\vert jj,ik>+(q-q^{-1})(q[2])^{1/2}\vert jk,ij>,$$

$$R\vert jj,ik>=(q-q^{-1})(q[2])^{1/2}\vert jk,ij>+\vert ik,jj>,$$

$$R\vert ik,jl>=(q-q^{-1})\vert kl,ij>+\vert jl,ik>,$$

$$R\vert ij,ij>=(q^{3}-q^{-1})\vert jj,ii>+q^{2}\vert ij,ij>,$$

$$R\vert ij,ii>=q^{2}\vert ii,ij>+(q^{4}-1)\vert ij,ii>,$$

$$R\vert jj,ij>=q^{2}\vert ij,jj>+(q^{4}-1)\vert jj,ij>,$$

$$R\vert il,jk>=\vert jk,ij> +(q-q^{-1})\vert jl,ik>+(q^{2}-1)\vert kl,ij>,$$

$$R\vert jk,il>=(q-q^{-1})\vert jl,ik>+\vert il,jk>+(q^{2}-1)\vert kl,ij>,$$

$$R\vert jl,ik>=\vert ik,jl>+(q-q^{-1})(\vert jk,il>+il,jk>$$

$$+(q-q^{-1})^{2}\vert jl,ik>+(q^{3}-q)\vert kl,ij>,$$

$$R\vert jj,ii>=(q^{4}-q^{2}-1+q^{-2})\vert jj,ii>+\vert ii,jj>+
(q^{3}-q^{-1})\vert ij,ij>,$$

$$R\vert jk,ij>=(2q^{2}+q^{-2}-3)\vert jk,ij>+q\vert ij,jk>$$

$$+(q-q^{-1})(q[2])^{1/2}(\vert ik,jj>+\vert jj,ik>,$$

$$R\vert kl,ij>=(q-q^{-1})^{2}(q^{2}+1)\vert kl,ij>+(q^{3}-q)\vert
jl,ik>+(q^{2}-1)\vert jk,il>$$

$$+(q^{2}-1)\vert il,jk>+(q-q^{-1})\vert ik,jl>+\vert ij,lk>,$$

$$R\vert ij,ik>=q\vert ik,ij>+(q^{3/2}-q^{-1/2})[2]\vert jk,ii>,$$

$$R\vert ik,ij>=q\vert ij,ik>+(q^{2}-1)\vert ik,ij>
+(q^{2}-1)(q[2])^{1/2}\vert jk,ii>,$$

$$R\vert jk,ii>=(q-q^{-1})^{2}[2]q\vert jk,ii>+\vert ii,jk>$$

$$+(q^{2}-1)(q[2])^{1/2}\vert ik,ij> +(q-q^{-1})\vert ij,ik>,$$

$$R\vert kk,ij>=q\vert ij,kk>+q^{3/2}(q-q^{1-})^{2}[2]^{2}\vert kk,ij>$$

$$+(q^{2}-1)(q[2])^{1/2}\vert jk,ik>+(q-q^{-1})(q[2])^{1/2}\vert ik,jk>,$$

$$R\vert ik,jk>=q\vert jk,ik>+(q-q^{-1})([2]q)^{1/2}\vert kk,ij>,$$

$$R\vert jk,ik>=(q^{2}-1)\vert jk,ik>+q\vert ik,jk>+(q^{2}-1)(q[2])^{1/2}\vert kk,ij>.\eqno(37)$$

\vskip .3truecm
\noindent (3) $[1^{2}]\times [1^{2}]$ irrep.
\vskip .3truecm
$$R\left\vert
\begin{array}{l}
i\\
j
\end{array}\right.,
\left.\begin{array}{l}
i\\
j
\end{array}\right> =q^{2}\left\vert\begin{array}{l}
i\\
j
\end{array}\right.,\left.\begin{array}{l}
i\\
j
\end{array}\right>,$$

$$R\left\vert\begin{array}{l}
i\\
j
\end{array}\right.,\left.\begin{array}{l}
i\\
k
\end{array}\right> =q\left\vert\begin{array}{l}
i\\
k
\end{array}\right.,\left.\begin{array}{l}
i\\
j
\end{array}\right>,$$

$$R\left\vert\begin{array}{l}
i\\
j
\end{array}\right.,\left.\begin{array}{l}
j\\
k
\end{array}\right>=q\left\vert\begin{array}{l}
j\\
k
\end{array}\right.,\left.\begin{array}{l}
i\\
j
\end{array}\right>,$$

$$R\left\vert\begin{array}{l}
i\\
k
\end{array}\right.,\left.\begin{array}{l}
j\\
k
\end{array}\right>=q\left\vert\begin{array}{l}
j\\
k
\end{array}\right.,\left.\begin{array}{l}
i\\
k
\end{array}\right>,$$

$$R\left\vert\begin{array}{l}
i\\
k
\end{array}\right.,\left.\begin{array}{l}
j\\
l
\end{array}\right>=\left\vert\begin{array}{l}
j\\
l
\end{array}\right.,\left.\begin{array}{l}
i\\
k
\end{array}\right>+(q-q^{-1})\left\vert\begin{array}{l}
k\\
l
\end{array}\right.,\left.\begin{array}{l}
i\\
j
\end{array}\right>,$$

$$R\left\vert\begin{array}{l}
i\\
k
\end{array}\right.,\left.\begin{array}{l}
i\\
j
\end{array}\right> =(q^{2}-1)\left\vert\begin{array}{l}
i\\
k
\end{array}\right.,\left.\begin{array}{l}
i\\
j
\end{array}\right>+q
\left\vert\begin{array}{l}
i\\
j
\end{array}\right.,\left.\begin{array}{l}
i\\
k
\end{array}\right>,$$

$$R\left\vert\begin{array}{l}
j\\
k
\end{array}\right.,\left.\begin{array}{l}
i\\
j
\end{array}\right>=q\left\vert\begin{array}{l}
i\\
j
\end{array}\right.,\left.\begin{array}{l}
j\\
k
\end{array}\right>+(q^{2}-1)\left\vert\begin{array}{l}
j\\
k
\end{array}\right.,\left.\begin{array}{l}
i\\
j
\end{array}\right>,$$

$$R\left\vert\begin{array}{l}
j\\
k
\end{array}\right.,\left.\begin{array}{l}
i\\
k
\end{array}\right>=q\left\vert\begin{array}{l}
i\\
k
\end{array}\right.,\left.\begin{array}{l}
j\\
k
\end{array}\right>+(q^{2}-1)\left\vert\begin{array}{l}
j\\
k
\end{array}\right.,\left.\begin{array}{l}
i\\
k
\end{array}\right>,$$

$$R\left\vert\begin{array}{l}
k\\
l
\end{array}\right.,\left.\begin{array}{l}
i\\
j
\end{array}\right>=q^{-1}[2](q-q^{-1})^{2}\left\vert\begin{array}{l}
k\\
l
\end{array}\right.,\left.\begin{array}{l}
i\\
j
\end{array}\right>+(q^{-1}-q^{-3})\left\vert\begin{array}{l}
j\\
l
\end{array}\right.,\left.\begin{array}{l}
i\\
k
\end{array}\right>-(1-q^{-2})\left\vert\begin{array}{l}
j\\
k
\end{array}\right.,\left.\begin{array}{l}
i\\
l
\end{array}\right>$$

$$-(1-q^{-2})\left\vert\begin{array}{l}
i\\
l
\end{array}\right.,\left.\begin{array}{l}
j\\
k
\end{array}\right>+(q-q^{-1})\left\vert\begin{array}{l}
i\\
k
\end{array}\right.,\left.\begin{array}{l}
j\\
l
\end{array}\right>+\left\vert\begin{array}{l}
i\\
j
\end{array}\right.,\left.\begin{array}{l}
k\\
l
\end{array}\right>,$$

$$R\left\vert\begin{array}{l}
j\\
l
\end{array}\right.,\left.\begin{array}{l}
i\\
k
\end{array}\right>=(q-q^{-1})^{2}\left\vert\begin{array}{l}
j\\
l
\end{array}\right.,\left.\begin{array}{l}
i\\
k
\end{array}\right>+(q^{-1}-q^{-3})\left\vert\begin{array}{l}
k\\
l
\end{array}\right.,\left.\begin{array}{l}
i\\
j
\end{array}\right>+(q-q^{-1})\left\vert\begin{array}{l}
j\\
k
\end{array}\right.,\left.\begin{array}{l}
i\\
l
\end{array}\right>$$

$$+(q-q^{-1})\left\vert\begin{array}{l}
i\\
l
\end{array}\right.,\left.\begin{array}{l}
j\\
k
\end{array}\right>+\left\vert\begin{array}{l}
i\\
k
\end{array}\right.,\left.\begin{array}{l}
j\\
l
\end{array}\right>,$$

$$R\left\vert\begin{array}{l}
i\\
l
\end{array}\right.,\left.\begin{array}{l}
j\\
k
\end{array}\right>=\left\vert\begin{array}{l}
j\\
k
\end{array}\right.,\left.\begin{array}{l}
i\\
l
\end{array}\right>+(q-q^{-1})\left\vert\begin{array}{l}
j\\
l
\end{array}\right.,\left.\begin{array}{l}
i\\
k
\end{array}\right>-(1-q^{-2})\left\vert\begin{array}{l}
k\\
l
\end{array}\right.,\left.\begin{array}{l}
i\\
j
\end{array}\right>,$$

$$R\left\vert\begin{array}{l}
j\\
k
\end{array}\right.,\left.\begin{array}{l}
i\\
l
\end{array}\right>=\left\vert\begin{array}{l}
i\\
l
\end{array}\right.,\left.\begin{array}{l}
j\\
k
\end{array}\right>-(1-q^{-2})\left\vert\begin{array}{l}
k\\
l
\end{array}\right.,\left.\begin{array}{l}
i\\
j
\end{array}\right>+(q-q^{-1})\left\vert\begin{array}{l}
j\\
l
\end{array}\right.,\left.\begin{array}{l}
i\\
k
\end{array}\right>.\eqno(38)$$
\vskip .3truecm
\noindent (4) $[21]\times [21]$ irrep.
\vskip .3truecm
In this case $R=g_{3}g_{4}g_{2}g_{3}g_{1}g_{2}g_{5}g_{4}g_{3}$.
Using the CGC's for $[1]\times [1]\downarrow [2]$ or $[11]$, and
$[2]\times [1]\downarrow [21]$ or $[11]\times [1]\downarrow [21]$, we
have the following expansions for uncoupled basis vectors of $[21]\times [21]$.
\vskip .3truecm
$$\left\vert\begin{array}{l}
ij\\
j
\end{array}\right>=\sqrt{{q\over{[3]!}}}(1+qg_{1}-q^{-1}[2]g_{2}g_{1})\vert i,j,j> =B^{0}_{12}\vert i,j,j>,$$

$$\left\vert\begin{array}{l}
ii\\
j
\end{array}\right>=\sqrt{{q\over{[3]!}}}([2]-q^{-2}g_{2}-q^{-1}g_{1}g_{2})\vert i,i,j> =B^{1}_{12}\vert i,j,j>,$$

$$\left\vert\begin{array}{l}
ik\\
j
\end{array}\right>={1\over{[2]}}(g_{2}+qg_{1}g_{2}-q^{-1}g_{2}g_{1}-g_{2}g_{1}g_{2})\vert i,j,k> =B^{2}_{12}\vert i,j,k>,$$

$$\left\vert\begin{array}{l}
ij\\
k
\end{array}\right>=\sqrt{{1\over{[3]}}}(1+qg_{1}-{q^{-2}\over{[2]}}g_{2}
-{q^{-1}\over{[2]}}g_{1}g_{2}
-{q^{-1}\over{[2]}}g_{2}g_{1}
-{1\over{[2]}}g_{1}g_{2}g_{1})\vert i,j,k> =B^{3}_{12}\vert i,j,k>.\eqno(39)$$
\vskip .3truecm
\noindent  We can also prove that
\vskip .3truecm
$${R}B_{12}^{p}B_{45}^{r}=B_{45}^{p}B_{12}^{r}R~~{\rm for\ }~~ 0\leq p,r\leq 3,\eqno(40a)$$
\vskip .3truecm
\noindent and
\vskip .3truecm
$${R}B_{12}^{p}=B_{45}^{p}R,~~{R}B_{45}^{p}=B^{p}_{12}R~~{\rm for \ } ~~0\leq p\leq 3.\eqno(40b)$$
\vskip .3truecm
After a long calculation with the help of Hecke algebra relations defined
by (17), We can also derive the following quadratic equation for R
\vskip .3truecm
$$R^{2}=(q-q^{-1})^{2}{R}g_{1}g_{2}g_{5}g_{4}g_{3}+(q-q^{-1})^{3}R+$$

$$(q-q^{-1})^{2}{R}g_{1}g_{2}g_{1}g_{5}g_{4}g_{5}g_{3}+(q-q^{-1})^{3}g_{2}{R}g_{2}+
(q-q^{-1})^{3}g_{5}{R}g_{4}g_{5}g_{1}+(q-q^{-1})^{3}{R}g_{1}g_{4}+$$

$$(q-q^{-1})^{2}(g_{2}g_{3}g_{4}g_{5}g_{4}g_{1}g_{2}g_{1}g_{3}g_{2}+
g_{2}g_{3}g_{4}g_{5}g_{4}g_{3}g_{2}g_{4}g_{3}g_{4}+
g_{2}g_{3}g_{4}g_{5}g_{4}g_{3}g_{2}g_{3}+$$

$$g_{3}g_{1}g_{2}g_{1}g_{3}g_{4}g_{5}g_{4}g_{3}g_{2}g_{3}g_{1}+
g_{3}g_{4}g_{5}g_{4}g_{1}g_{3}g_{2}g_{3}g_{4}g_{3}g_{5}g_{1}+
g_{3}g_{4}g_{3}g_{5}g_{1}g_{2}g_{1}g_{4}g_{3}g_{4}+$$

$$g_{3}g_{4}g_{1}g_{2}g_{3}g_{1}g_{2}g_{4}g_{3}g_{1}+
g_{3}g_{4}g_{2}g_{1}g_{5}g_{2}g_{3}+
g_{3}g_{4}g_{2}g_{3}g_{5}g_{2}g_{4}g_{3}+g_{3}g_{4}g_{2}g_{5}g_{4}g_{3}+$$

$$g_{3}g_{4}g_{2}g_{3}g_{1}g_{2}g_{4}g_{3}+g_{3}g_{2}g_{4}g_{3}+
g_{5}g_{3}g_{4}g_{1}g_{2}g_{3}g_{4}g_{5}g_{2}g_{1}+g_{4}g_{3}g_{2}g_{4}g_{3}g_{4})$$

$$+(q-q^{-1})(g_{2}g_{3}g_{4}g_{5}g_{4}g_{3}g_{2}+g_{3}g_{4}g_{1}g_{2}g_{3}g_{2}
g_{4}g_{3}g_{1}+g_{3}g_{4}g_{5}g_{4}g_{3}+$$

$$g_{3}g_{4}g_{1}g_{2}g_{1}g_{4}g_{3}+g_{3}g_{4}g_{3}+g_{3}+g_{3}g_{2}g_{3}
+g_{5}g_{4}g_{1}g_{2}g_{3}g_{4}g_{5}g_{2}g_{1}$$

$$+g_{3}g_{4}g_{2}g_{3}g_{4})+1.\eqno(41)$$
\vskip .3truecm
\noindent Using these relations, we  get
\vskip .3truecm
$$R\left\vert
\begin{array}{l}
ii\\
j
\end{array}\right.,\left.\begin{array}{l}
ii\\
j
\end{array}\right>=q^{5}\left\vert
\begin{array}{l}
ii\\
j
\end{array}\right.,\left.\begin{array}{l}
ii\\
j
\end{array}\right>$$

$$R\left\vert
\begin{array}{l}
ij\\
j
\end{array}\right.,\left.\begin{array}{l}
ij\\
j
\end{array}\right>=q^{5}\left\vert
\begin{array}{l}
ij\\
j
\end{array}\right.,\left.\begin{array}{l}
ij\\
j
\end{array}\right>$$

$$R\left\vert
\begin{array}{l}
ii\\
j
\end{array}\right.,\left.\begin{array}{l}
ij\\
j
\end{array}\right>=q^{4}\left\vert
\begin{array}{l}
ij\\
j
\end{array}\right.,\left.\begin{array}{l}
ii\\
j
\end{array}\right>$$

$$R\left\vert
\begin{array}{l}
ik\\
k
\end{array}\right.,\left.\begin{array}{l}
jk\\
k
\end{array}\right>=q^{4}\left\vert
\begin{array}{l}
jk\\
k
\end{array}\right.,\left.\begin{array}{l}
ik\\
k
\end{array}\right>$$

$$R\left\vert
\begin{array}{l}
ij\\
j
\end{array}\right.,\left.\begin{array}{l}
jj\\
k
\end{array}\right>=q^{4}\left\vert
\begin{array}{l}
jj\\
k
\end{array}\right.,\left.\begin{array}{l}
ij\\
j
\end{array}\right>$$

$$R\left\vert
\begin{array}{l}
ii\\
j
\end{array}\right.,\left.\begin{array}{l}
ii\\
k
\end{array}\right>=q^{4}\left\vert
\begin{array}{l}
ii\\
k
\end{array}\right.,\left.\begin{array}{l}
ii\\
j
\end{array}\right>$$

$$R\left\vert
\begin{array}{l}
ij\\
j
\end{array}\right.,\left.\begin{array}{l}
ik\\
j
\end{array}\right>=q^{3}\left\vert
\begin{array}{l}
ik\\
j
\end{array}\right.,\left.\begin{array}{l}
ij\\
j
\end{array}\right>$$

$$R\left\vert
\begin{array}{l}
ik\\
j
\end{array}\right.,\left.\begin{array}{l}
ik\\
k
\end{array}\right>=q^{3}\left\vert
\begin{array}{l}
ik\\
k
\end{array}\right.,\left.\begin{array}{l}
ik\\
j
\end{array}\right>$$

$$R\left\vert
\begin{array}{l}
ii\\
j
\end{array}\right.,\left.\begin{array}{l}
ij\\
k
\end{array}\right>=q^{3}\left\vert
\begin{array}{l}
ij\\
k
\end{array}\right.,\left.\begin{array}{l}
ii\\
j
\end{array}\right>$$

$$R\left\vert
\begin{array}{l}
ij\\
k
\end{array}\right.,\left.\begin{array}{l}
jk\\
k
\end{array}\right>=q^{3}\left\vert
\begin{array}{l}
jk\\
k
\end{array}\right.,\left.\begin{array}{l}
ij\\
k
\end{array}\right>$$

$$R\left\vert
\begin{array}{l}
ii\\
j
\end{array}\right.,\left.\begin{array}{l}
ik\\
j
\end{array}\right>=q^{3}\left\vert
\begin{array}{l}
ik\\
j
\end{array}\right.,\left.\begin{array}{l}
ii\\
j
\end{array}\right>$$

$$R\left\vert
\begin{array}{l}
ik\\
j
\end{array}\right.,\left.\begin{array}{l}
jk\\
k
\end{array}\right>=q^{3}\left\vert
\begin{array}{l}
jk\\
k
\end{array}\right.,\left.\begin{array}{l}
ik\\
j
\end{array}\right>$$

$$R\left\vert
\begin{array}{l}
ii\\
j
\end{array}\right.,\left.\begin{array}{l}
jj\\
k
\end{array}\right>=q^{2}\left\vert
\begin{array}{l}
jj\\
k
\end{array}\right.,\left.\begin{array}{l}
ii\\
j
\end{array}\right>$$

$$R\left\vert
\begin{array}{l}
ii\\
k
\end{array}\right.,\left.\begin{array}{l}
jk\\
k
\end{array}\right>=q^{2}\left\vert
\begin{array}{l}
jk\\
k
\end{array}\right.,\left.\begin{array}{l}
ii\\
k
\end{array}\right>$$

$$R\left\vert
\begin{array}{l}
ii\\
j
\end{array}\right.,\left.\begin{array}{l}
ik\\
k
\end{array}\right>=q^{2}\left\vert
\begin{array}{l}
ik\\
k
\end{array}\right.,\left.\begin{array}{l}
ii\\
j
\end{array}\right>$$

$$R\left\vert
\begin{array}{l}
ij\\
j
\end{array}\right.,\left.\begin{array}{l}
jk\\
k
\end{array}\right>=q^{2}\left\vert
\begin{array}{l}
jk\\
k
\end{array}\right.,\left.\begin{array}{l}
ij\\
j
\end{array}\right>$$

$$R\left\vert
\begin{array}{l}
ii\\
j
\end{array}\right.,\left.\begin{array}{l}
jk\\
k
\end{array}\right>=q\left\vert
\begin{array}{l}
jk\\
k
\end{array}\right.,\left.\begin{array}{l}
ii\\
j
\end{array}\right>$$

$$R\left\vert
\begin{array}{l}
ij\\
j
\end{array}\right.,\left.\begin{array}{l}
ik\\
k
\end{array}\right>=q\left\vert
\begin{array}{l}
ik\\
k
\end{array}\right.,\left.\begin{array}{l}
ij\\
j
\end{array}\right>+(q^{2}-1)\left\vert
\begin{array}{l}
jk\\
k
\end{array}\right.,\left.\begin{array}{l}
ii\\
j
\end{array}\right>$$

$$R\left\vert
\begin{array}{l}
ii\\
k
\end{array}\right.,\left.\begin{array}{l}
jj\\
k
\end{array}\right>=q\left\vert
\begin{array}{l}
jj\\
k
\end{array}\right.,\left.\begin{array}{l}
ii\\
k
\end{array}\right>+(q^{2}-1)\left\vert
\begin{array}{l}
jk\\
k
\end{array}\right.,\left.\begin{array}{l}
ii\\
j
\end{array}\right>$$

$$R\left\vert
\begin{array}{l}
ij\\
k
\end{array}\right.,\left.\begin{array}{l}
jj\\
k
\end{array}\right>=q^{3}\left\vert
\begin{array}{l}
jj\\
k
\end{array}\right.,\left.\begin{array}{l}
ij\\
k
\end{array}\right>+(q^{7/2}-q^{3/2})[2]^{-1/2}\left\vert
\begin{array}{l}
jk\\
k
\end{array}\right.,\left.\begin{array}{l}
ij\\
j
\end{array}\right>$$

$$R\left\vert
\begin{array}{l}
ii\\
k
\end{array}\right.,\left.\begin{array}{l}
ij\\
k
\end{array}\right>=q^{3}\left\vert
\begin{array}{l}
ij\\
k
\end{array}\right.,\left.\begin{array}{l}
ii\\
k
\end{array}\right>+(q^{7/2}-q^{3/2})[2]^{-1/2}\left\vert
\begin{array}{l}
ik\\
k
\end{array}\right.,\left.\begin{array}{l}
ii\\
j
\end{array}\right>$$

$$R\left\vert
\begin{array}{l}
ii\\
k
\end{array}\right.,\left.\begin{array}{l}
ik\\
j
\end{array}\right>=q^{3}\left\vert
\begin{array}{l}
ik\\
j
\end{array}\right.,\left.\begin{array}{l}
ii\\
k
\end{array}\right>+(q^{7/2}-q^{3/2})\left({[3]\over{[2]}}\right)^{1/2}\left\vert
\begin{array}{l}
ik\\
k
\end{array}\right.,\left.\begin{array}{l}
ii\\
j
\end{array}\right>$$

$$R\left\vert
\begin{array}{l}
ij\\
k
\end{array}\right.,\left.\begin{array}{l}
jj\\
k
\end{array}\right>=q^{3}\left\vert
\begin{array}{l}
jj\\
k
\end{array}\right.,\left.\begin{array}{l}
ik\\
j
\end{array}\right>+(q^{7/2}-q^{3/2})\left({[3]\over{[2]}}\right)^{1/2}\left\vert
\begin{array}{l}
jk\\
k
\end{array}\right.,\left.\begin{array}{l}
ij\\
j
\end{array}\right>$$

$$R\left\vert
\begin{array}{l}
ij\\
j
\end{array}\right.,\left.\begin{array}{l}
ii\\
j
\end{array}\right>=(q^{5}-q^{3})\left\vert
\begin{array}{l}
ij\\
j
\end{array}\right.,\left.\begin{array}{l}
ii\\
j
\end{array}\right>+q^{4}\left\vert
\begin{array}{l}
ii\\
j
\end{array}\right.,\left.\begin{array}{l}
ij\\
j
\end{array}\right>$$

$$R\left\vert
\begin{array}{l}
ii\\
k
\end{array}\right.,\left.\begin{array}{l}
ii\\
j
\end{array}\right>=(q^{5}-q^{3})\left\vert
\begin{array}{l}
ii\\
k
\end{array}\right.,\left.\begin{array}{l}
ii\\
j
\end{array}\right>+q^{4}\left\vert
\begin{array}{l}
ii\\
j
\end{array}\right.,\left.\begin{array}{l}
ii\\
k
\end{array}\right>$$

$$R\left\vert
\begin{array}{l}
jj\\
k
\end{array}\right.,\left.\begin{array}{l}
ij\\
j
\end{array}\right>=(q^{5}-q^{3})\left\vert
\begin{array}{l}
jj\\
k
\end{array}\right.,\left.\begin{array}{l}
ij\\
j
\end{array}\right>+q^{4}\left\vert
\begin{array}{l}
ij\\
j
\end{array}\right.,\left.\begin{array}{l}
jj\\
k
\end{array}\right>$$

$$R\left\vert
\begin{array}{l}
jk\\
k
\end{array}\right.,\left.\begin{array}{l}
ik\\
k
\end{array}\right>=q^{4}\left\vert
\begin{array}{l}
ik\\
k
\end{array}\right.,\left.\begin{array}{l}
jk\\
k
\end{array}\right>+(q^{5}-q^{3})\left\vert
\begin{array}{l}
jk\\
k
\end{array}\right.,\left.\begin{array}{l}
ik\\
k
\end{array}\right>$$

$$R\left\vert
\begin{array}{l}
ij\\
j
\end{array}\right.,\left.\begin{array}{l}
ij\\
k
\end{array}\right>=q^{3}\left\vert
\begin{array}{l}
ij\\
k
\end{array}\right.,\left.\begin{array}{l}
ij\\
j
\end{array}\right>+(q^{7/2}-q^{3/2})[2]^{1/2}\left\vert
\begin{array}{l}
jj\\
k
\end{array}\right.,\left.\begin{array}{l}
ii\\
j
\end{array}\right>$$

$$R\left\vert
\begin{array}{l}
ij\\
k
\end{array}\right.,\left.\begin{array}{l}
ik\\
k
\end{array}\right>=q^{3}\left\vert
\begin{array}{l}
ik\\
k
\end{array}\right.,\left.\begin{array}{l}
ij\\
k
\end{array}\right>+(q^{7/2}-q^{3/2})[2]^{1/2}\left\vert
\begin{array}{l}
jk\\
k
\end{array}\right.,\left.\begin{array}{l}
ii\\
k
\end{array}\right>$$

$$R\left\vert
\begin{array}{l}
ij\\
k
\end{array}\right.,\left.\begin{array}{l}
ij\\
j
\end{array}\right>=q^{3}\left\vert
\begin{array}{l}
ij\\
j
\end{array}\right.,\left.\begin{array}{l}
ij\\
k
\end{array}\right>+(q^{4}-q^{2})[2]^{-1}\left\vert
\begin{array}{l}
ij\\
k
\end{array}\right.,\left.\begin{array}{l}
ij\\
j
\end{array}\right>+$$

$$(q^{4}-q^{2}){[3]^{1/2}\over{[2]}}\left\vert
\begin{array}{l}
ik\\
j
\end{array}\right.,\left.\begin{array}{l}
ij\\
j
\end{array}\right>+(q^{9/2}-q^{5/2})[2]^{-1/2}\left\vert
\begin{array}{l}
jj\\
k
\end{array}\right.,\left.\begin{array}{l}
ii\\
j
\end{array}\right>$$

$$R\left\vert
\begin{array}{l}
ik\\
k
\end{array}\right.,\left.\begin{array}{l}
ij\\
k
\end{array}\right>=q^{3}\left\vert
\begin{array}{l}
ij\\
k
\end{array}\right.,\left.\begin{array}{l}
ik\\
k
\end{array}\right>+(q^{4}-q^{2})[2]^{-1}\left\vert
\begin{array}{l}
ik\\
k
\end{array}\right.,\left.\begin{array}{l}
ij\\
k
\end{array}\right>+$$

$$(q^{4}-q^{2}){[3]^{1/2}\over{[2]}}\left\vert
\begin{array}{l}
ik\\
k
\end{array}\right.,\left.\begin{array}{l}
ik\\
j
\end{array}\right>+(q^{9/2}-q^{5/2})[2]^{-1/2}\left\vert
\begin{array}{l}
jk\\
k
\end{array}\right.,\left.\begin{array}{l}
ii\\
k
\end{array}\right>$$

$$R\left\vert
\begin{array}{l}
ij\\
j
\end{array}\right.,\left.\begin{array}{l}
ii\\
k
\end{array}\right>=q^{2}\left\vert
\begin{array}{l}
ii\\
k
\end{array}\right.,\left.\begin{array}{l}
ij\\
j
\end{array}\right>+(q^{7/2}-q^{3/2})[2]^{1/2}\left\vert
\begin{array}{l}
ij\\
k
\end{array}\right.,\left.\begin{array}{l}
ii\\
j
\end{array}\right>$$

$$R\left\vert
\begin{array}{l}
jj\\
k
\end{array}\right.,\left.\begin{array}{l}
ik\\
k
\end{array}\right>=q^{2}\left\vert
\begin{array}{l}
ik\\
k
\end{array}\right.,\left.\begin{array}{l}
jj\\
k
\end{array}\right>+(q^{7/2}-q^{3/2})[2]^{1/2}\left\vert
\begin{array}{l}
jk\\
k
\end{array}\right.,\left.\begin{array}{l}
ij\\
k
\end{array}\right>$$

$$R\left\vert
\begin{array}{l}
ii\\
k
\end{array}\right.,\left.\begin{array}{l}
ij\\
j
\end{array}\right>=q^{2}\left\vert
\begin{array}{l}
ij\\
j
\end{array}\right.,\left.\begin{array}{l}
ii\\
k
\end{array}\right>+(q^{7/2}-q^{3/2})[2]^{-1/2}\left\vert
\begin{array}{l}
ij\\
k
\end{array}\right.,\left.\begin{array}{l}
ii\\
j
\end{array}\right>+$$

$$(q^{7/2}-q^{3/2})\left({[3]\over{[2]}}\right)^{1/2}\left\vert
\begin{array}{l}
ik\\
j
\end{array}\right.,\left.\begin{array}{l}
ii\\
j
\end{array}\right>$$

$$R\left\vert
\begin{array}{l}
ik\\
k
\end{array}\right.,\left.\begin{array}{l}
jj\\
k
\end{array}\right>=q^{2}\left\vert
\begin{array}{l}
jj\\
k
\end{array}\right.,\left.\begin{array}{l}
ik\\
k
\end{array}\right>+(q^{7/2}-q^{3/2})[2]^{-1/2}\left\vert
\begin{array}{l}
jk\\
k
\end{array}\right.,\left.\begin{array}{l}
ij\\
k
\end{array}\right>+$$

$$(q^{7/2}-q^{3/2})\left({[3]\over{[2]}}\right)^{1/2}\left\vert
\begin{array}{l}
jk\\
k
\end{array}\right.,\left.\begin{array}{l}
ik\\
j
\end{array}\right>$$

$$R\left\vert
\begin{array}{l}
ij\\
k
\end{array}\right.,\left.\begin{array}{l}
ii\\
j
\end{array}\right>=q^{3}\left\vert
\begin{array}{l}
ii\\
j
\end{array}\right.,\left.\begin{array}{l}
ij\\
k
\end{array}\right>+(q^{7/2}-q^{3/2})[2]^{1/2}\left\vert
\begin{array}{l}
ij\\
j
\end{array}\right.,\left.\begin{array}{l}
ii\\
k
\end{array}\right>+$$

$$(q^{7/2}-q^{3/2})[2]^{-1/2}\left\vert
\begin{array}{l}
ii\\
k
\end{array}\right.,\left.\begin{array}{l}
ij\\
j
\end{array}\right>+(1-2q^{2}+q^{6})[2]^{-1}\left\vert
\begin{array}{l}
ij\\
k
\end{array}\right.,\left.\begin{array}{l}
ii\\
j
\end{array}\right>+(1-q^{2}){[3]^{1/2}\over{[2]}}\left\vert
\begin{array}{l}
ik\\
j
\end{array}\right.,\left.\begin{array}{l}
ii\\
j
\end{array}\right>$$

$$R\left\vert
\begin{array}{l}
jk\\
k
\end{array}\right.,\left.\begin{array}{l}
ij\\
k
\end{array}\right>=q^{3}\left\vert
\begin{array}{l}
ij\\
k
\end{array}\right.,\left.\begin{array}{l}
jk\\
k
\end{array}\right>+(q^{7/2}-q^{3/2})[2]^{1/2}\left\vert
\begin{array}{l}
ij\\
k
\end{array}\right.,\left.\begin{array}{l}
ik\\
k
\end{array}\right>+$$

$$(q^{7/2}-q^{3/2})[2]^{-1/2}\left\vert
\begin{array}{l}
ik\\
k
\end{array}\right.,\left.\begin{array}{l}
jj\\
k
\end{array}\right>+(1-2q^{2}+q^{6})[2]^{-1}\left\vert
\begin{array}{l}
jk\\
k
\end{array}\right.,\left.\begin{array}{l}
ij\\
k
\end{array}\right>+(1-q^{2}){[3]^{1/2}\over{[2]}}\left\vert
\begin{array}{l}
jk\\
k
\end{array}\right.,\left.\begin{array}{l}
ik\\
j
\end{array}\right>$$

$$R\left\vert
\begin{array}{l}
ik\\
j
\end{array}\right.,\left.\begin{array}{l}
ii\\
j
\end{array}\right>=q^{3}\left\vert
\begin{array}{l}
ii\\
j
\end{array}\right.,\left.\begin{array}{l}
ik\\
j
\end{array}\right>+(q^{7/2}-q^{3/2})\left({[3]\over{[2]}}\right)^{1/2}\left\vert
\begin{array}{l}
ii\\
k
\end{array}\right.,\left.\begin{array}{l}
ij\\
j
\end{array}\right>+$$

$$(1-q^{2}){[3]^{1/2}\over{[2]}}\left\vert
\begin{array}{l}
ij\\
k
\end{array}\right.,\left.\begin{array}{l}
ii\\
j
\end{array}\right>+(q^{6}-1)[2]^{-1}\left\vert
\begin{array}{l}
ik\\
j
\end{array}\right.,\left.\begin{array}{l}
ii\\
j
\end{array}\right>$$

$$R\left\vert
\begin{array}{l}
jk\\
k
\end{array}\right.,\left.\begin{array}{l}
ik\\
j
\end{array}\right>=q^{3}\left\vert
\begin{array}{l}
ik\\
j
\end{array}\right.,\left.\begin{array}{l}
jk\\
k
\end{array}\right>+(q^{7/2}-q^{3/2})\left({[3]^{1/2}\over{[2]}}\right)^{1/2}\left\vert
\begin{array}{l}
ik\\
k
\end{array}\right.,\left.\begin{array}{l}
jj\\
k
\end{array}\right>+$$

$$(1-q^{2}){[3]^{1/2}\over{[2]}}\left\vert
\begin{array}{l}
jk\\
k
\end{array}\right.,\left.\begin{array}{l}
ij\\
k
\end{array}\right>+(q^{6}-1)[2]^{-1}\left\vert
\begin{array}{l}
jk\\
k
\end{array}\right.,\left.\begin{array}{l}
ik\\
j
\end{array}\right>$$

$$R\left\vert
\begin{array}{l}
ik\\
j
\end{array}\right.,\left.\begin{array}{l}
ij\\
j
\end{array}\right>=q^{3}\left\vert
\begin{array}{l}
ij\\
j
\end{array}\right.,\left.\begin{array}{l}
ik\\
j
\end{array}\right>+(q^{4}-q^{2}){[3]^{1/2}\over{[2]}}\left\vert
\begin{array}{l}
ij\\
k
\end{array}\right.,\left.\begin{array}{l}
ij\\
j
\end{array}\right>+(q^{6}-1)[2]^{-1}\left\vert
\begin{array}{l}
ik\\
j
\end{array}\right.,\left.\begin{array}{l}
ij\\
j
\end{array}\right>+$$

$$(q^{1/2}-q^{5/2})\left({[3]\over{[2]}}\right)^{1/2}\left\vert
\begin{array}{l}
jj\\
k
\end{array}\right.,\left.\begin{array}{l}
ii\\
j
\end{array}\right>$$

$$R\left\vert
\begin{array}{l}
ik\\
k
\end{array}\right.,\left.\begin{array}{l}
ik\\
j
\end{array}\right>=q^{3}\left\vert
\begin{array}{l}
ik\\
j
\end{array}\right.,\left.\begin{array}{l}
ik\\
k
\end{array}\right>+(q^{4}-q^{2}){[3]^{1/2}\over{[2]}}\left\vert
\begin{array}{l}
ik\\
k
\end{array}\right.,\left.\begin{array}{l}
ij\\
k
\end{array}\right>+$$

$$(q^{6}-1)[2]^{-1}\left\vert
\begin{array}{l}
ik\\
k
\end{array}\right.,\left.\begin{array}{l}
ik\\
j
\end{array}\right>+(q^{1/2}-q^{5/2})\left({[3]\over{[2]}}\right)^{1/2}\left\vert
\begin{array}{l}
jk\\
k
\end{array}\right.,\left.\begin{array}{l}
ii\\
k
\end{array}\right>$$

$$R\left\vert
\begin{array}{l}
jj\\
k
\end{array}\right.,\left.\begin{array}{l}
ii\\
j
\end{array}\right>=q^{2}\left\vert
\begin{array}{l}
ii\\
j
\end{array}\right.,\left.\begin{array}{l}
jj\\
k
\end{array}\right>+(q^{7/2}-q^{3/2})[2]^{1/2}\left\vert
\begin{array}{l}
ij\\
j
\end{array}\right.,\left.\begin{array}{l}
ij\\
k
\end{array}\right>+$$

$$(q^{9/2}-q^{5/2})[2]^{-1/2}\left\vert
\begin{array}{l}
ij\\
k
\end{array}\right.,\left.\begin{array}{l}
ij\\
j
\end{array}\right>+(q^{1/2}-q^{5/2})\left({[3]\over{[2]}}\right)^{1/2}\left\vert
\begin{array}{l}
ik\\
j
\end{array}\right.,\left.\begin{array}{l}
ij\\
j
\end{array}\right>+(q-2q^{3}+q^{5})\left\vert
\begin{array}{l}
jj\\
k
\end{array}\right.,\left.\begin{array}{l}
ii\\
j
\end{array}\right>$$

$$R\left\vert
\begin{array}{l}
jk\\
k
\end{array}\right.,\left.\begin{array}{l}
ii\\
k
\end{array}\right>=q^{2}\left\vert
\begin{array}{l}
ii\\
k
\end{array}\right.,\left.\begin{array}{l}
jk\\
k
\end{array}\right>+(q^{7/2}-q^{3/2})[2]^{1/2}\left\vert
\begin{array}{l}
ij\\
k
\end{array}\right.,\left.\begin{array}{l}
ik\\
k
\end{array}\right>+$$

$$(q^{9/2}-q^{5/2})[2]^{-1/2}\left\vert
\begin{array}{l}
ik\\
k
\end{array}\right.,\left.\begin{array}{l}
ij\\
k
\end{array}\right>+(q^{1/2}-q^{5/2})\left({[3]\over{[2]}}\right)^{1/2}\left\vert
\begin{array}{l}
ik\\
k
\end{array}\right.,\left.\begin{array}{l}
ik\\
j
\end{array}\right>+(q-2q^{3}+q^{5})\left\vert
\begin{array}{l}
jk\\
k
\end{array}\right.,\left.\begin{array}{l}
ii\\
k
\end{array}\right>$$

$$R\left\vert
\begin{array}{l}
ij\\
k
\end{array}\right.,\left.\begin{array}{l}
ii\\
k
\end{array}\right>=q^{3}\left\vert
\begin{array}{l}
ii\\
k
\end{array}\right.,\left.\begin{array}{l}
ij\\
k
\end{array}\right>+(q^{5}-q)\left\vert
\begin{array}{l}
ij\\
k
\end{array}\right.,\left.\begin{array}{l}
ii\\
k
\end{array}\right>+(q^{-1/2}-q^{3/2})[2]^{-1/2}\left\vert
\begin{array}{l}
ik\\
k
\end{array}\right.,\left.\begin{array}{l}
ii\\
j
\end{array}\right>$$

$$R\left\vert
\begin{array}{l}
jj\\
k
\end{array}\right.,\left.\begin{array}{l}
ij\\
k
\end{array}\right>=q^{3}\left\vert
\begin{array}{l}
ij\\
k
\end{array}\right.,\left.\begin{array}{l}
jj\\
k
\end{array}\right>+(q^{5}-q)\left\vert
\begin{array}{l}
jj\\
k
\end{array}\right.,\left.\begin{array}{l}
ij\\
k
\end{array}\right>+(q^{-1/2}-q^{3/2})[2]^{-1/2}\left\vert
\begin{array}{l}
jk\\
k
\end{array}\right.,\left.\begin{array}{l}
ij\\
j
\end{array}\right>$$

$$R\left\vert
\begin{array}{l}
ik\\
j
\end{array}\right.,\left.\begin{array}{l}
ii\\
k
\end{array}\right>=q^{3}\left\vert
\begin{array}{l}
ii\\
k
\end{array}\right.,\left.\begin{array}{l}
ik\\
j
\end{array}\right>+(q^{7/2}-q^{3/2})\left({[3]\over{[2]}}\right)^{-1/2}\left\vert
\begin{array}{l}
ik\\
k
\end{array}\right.,\left.\begin{array}{l}
ii\\
j
\end{array}\right>$$

$$R\left\vert
\begin{array}{l}
jj\\
k
\end{array}\right.,\left.\begin{array}{l}
ik\\
j
\end{array}\right>=q^{3}\left\vert
\begin{array}{l}
ik\\
j
\end{array}\right.,\left.\begin{array}{l}
jj\\
k
\end{array}\right>+(q^{7/2}-q^{3/2})\left({[3]\over{[2]}}\right)^{-1/2}\left\vert
\begin{array}{l}
jk\\
k
\end{array}\right.,\left.\begin{array}{l}
ij\\
j
\end{array}\right>$$

$$R\left\vert
\begin{array}{l}
ik\\
k
\end{array}\right.,\left.\begin{array}{l}
ii\\
j
\end{array}\right>=q^{2}\left\vert
\begin{array}{l}
ii\\
j
\end{array}\right.,\left.\begin{array}{l}
ik\\
k
\end{array}\right>+(q^{7/2}-q^{3/2}){[2]}^{-1/2}\left\vert
\begin{array}{l}
ii\\
k
\end{array}\right.,\left.\begin{array}{l}
ij\\
k
\end{array}\right>+ $$

$$(q^{-1/2}-q^{3/2}){[2]}^{-1/2}\left\vert
\begin{array}{l}
ij\\
k
\end{array}\right.,\left.\begin{array}{l}
ii\\
k
\end{array}\right> +(q^{7/2}-q^{3/2})\left({[3]\over{[2]}}\right)^{-1/2}\left\vert
\begin{array}{l}
ii\\
k
\end{array}\right.,\left.\begin{array}{l}
ik\\
j
\end{array}\right> +$$

$$(q^{7/2}-q^{3/2})\left({[3]\over{[2]}}\right)^{-1/2}\left\vert
\begin{array}{l}
ik\\
j
\end{array}\right.,\left.\begin{array}{l}
ik\\
k
\end{array}\right>+(q-2q^{3}+q^{5})\left\vert
\begin{array}{l}
ik\\
k
\end{array}\right.,\left.\begin{array}{l}
ii\\
j
\end{array}\right>$$

$$R\left\vert
\begin{array}{l}
jk\\
k
\end{array}\right.,\left.\begin{array}{l}
ij\\
j
\end{array}\right>=q^{2}\left\vert
\begin{array}{l}
ij\\
j
\end{array}\right.,\left.\begin{array}{l}
jk\\
k
\end{array}\right>+(q^{7/2}-q^{3/2}){[2]}^{-1/2}\left\vert
\begin{array}{l}
ij\\
k
\end{array}\right.,\left.\begin{array}{l}
jj\\
k
\end{array}\right>+$$

$$(q^{-1/2}-q^{3/2}){[2]}^{-1/2}\left\vert
\begin{array}{l}
jj\\
k
\end{array}\right.,\left.\begin{array}{l}
ij\\
k
\end{array}\right>$$

$$+(q^{7/2}-q^{3/2})\left({[3]\over{[2]}}\right)^{-1/2}\left\vert
\begin{array}{l}
ik\\
j
\end{array}\right.,\left.\begin{array}{l}
jj\\
k
\end{array}\right> +(q^{7/2}-q^{3/2})\left({[3]\over{[2]}}\right)^{-1/2}\left\vert
\begin{array}{l}
jj\\
k
\end{array}\right.,\left.\begin{array}{l}
ik\\
j
\end{array}\right>+(q-2q^{3}+q^{5})\left\vert
\begin{array}{l}
jk\\
k
\end{array}\right.,\left.\begin{array}{l}
ij\\
j
\end{array}\right>$$

$$R\left\vert
\begin{array}{l}
ij\\
k
\end{array}\right.,\left.\begin{array}{l}
ij\\
k
\end{array}\right>=q^{3}\left\vert
\begin{array}{l}
ij\\
k
\end{array}\right.,\left.\begin{array}{l}
ij\\
k
\end{array}\right>+(q^{-2}-2q^{2}+q^{4})[2]^{-1}\left\vert
\begin{array}{l}
jk\\
k
\end{array}\right.,\left.\begin{array}{l}
ii\\
j
\end{array}\right>$$

$$+(q^{4}-1)\left\vert
\begin{array}{l}
jj\\
k
\end{array}\right.,\left.\begin{array}{l}
ii\\
k
\end{array}\right>+(q^{3}-q)[2]^{-1}\left\vert
\begin{array}{l}
ik\\
k
\end{array}\right.,\left.\begin{array}{l}
ij\\
j
\end{array}\right>$$

$$R\left\vert
\begin{array}{l}
jk\\
k
\end{array}\right.,\left.\begin{array}{l}
ii\\
j
\end{array}\right>=q\left\vert
\begin{array}{l}
ii\\
j
\end{array}\right.,\left.\begin{array}{l}
jk\\
k
\end{array}\right>+(q^{2}-1)\left\vert
\begin{array}{l}
ij\\
j
\end{array}\right.,\left.\begin{array}{l}
ik\\
k
\end{array}\right>+(q^{-2}-2q^{2}+q^{4})[2]^{-1}\left\vert
\begin{array}{l}
ij\\
k
\end{array}\right.,\left.\begin{array}{l}
ij\\
k
\end{array}\right>+$$

$$(q^{4}-q^{2}){[3]^{1/2}\over{[2]}}\left\vert
\begin{array}{l}
ij\\
k
\end{array}\right.,\left.\begin{array}{l}
ik\\
j
\end{array}\right>+(q^{4}-q^{2}){[3]^{1/2}\over{[2]}}\left\vert
\begin{array}{l}
ik\\
j
\end{array}\right.,\left.\begin{array}{l}
ij\\
k
\end{array}\right>-(q^{4}-q^{-2})[2]^{-1}\left\vert
\begin{array}{l}
ik\\
j
\end{array}\right.,\left.\begin{array}{l}
ik\\
j
\end{array}\right>-$$

$$[2]^{2}\left(\left\vert
\begin{array}{l}
jj\\
k
\end{array}\right.,\left.\begin{array}{l}
ii\\
k
\end{array}\right>+\left\vert
\begin{array}{l}
ik\\
k
\end{array}\right.,\left.\begin{array}{l}
ij\\
j
\end{array}\right>\right)+(q^{-3}-2q^{-1}+2q-2q^{3}+q^{5})\left\vert
\begin{array}{l}
jk\\
k
\end{array}\right.,\left.\begin{array}{l}
ij\\
j
\end{array}\right>$$

$$R\left\vert
\begin{array}{l}
ik\\
j
\end{array}\right.,\left.\begin{array}{l}
ij\\
k
\end{array}\right>=q^{3}\left\vert
\begin{array}{l}
ij\\
k
\end{array}\right.,\left.\begin{array}{l}
ik\\
j
\end{array}\right>+(q^{3}-q){[3]^{1/2}\over{[2]}}\left\vert
\begin{array}{l}
ik\\
k
\end{array}\right.,\left.\begin{array}{l}
ij\\
j
\end{array}\right> +(q^{4}-q^{2}){[3]^{1/2}\over{[2]}}\left\vert
\begin{array}{l}
jk\\
k
\end{array}\right.,\left.\begin{array}{l}
ii\\
j
\end{array}\right>$$

$$R\left\vert
\begin{array}{l}
ij\\
k
\end{array}\right.,\left.\begin{array}{l}
ik\\
j
\end{array}\right>=q^{3}\left\vert
\begin{array}{l}
ik\\
j
\end{array}\right.,\left.\begin{array}{l}
ij\\
k
\end{array}\right>+(q^{3}-q){[3]^{1/2}\over{[2]}}\left\vert
\begin{array}{l}
ik\\
k
\end{array}\right.,\left.\begin{array}{l}
ij\\
j
\end{array}\right> +(q^{4}-q^{2}){[3]^{1/2}\over{[2]}}\left\vert
\begin{array}{l}
jk\\
k
\end{array}\right.,\left.\begin{array}{l}
ii\\
j
\end{array}\right>$$

$$R\left\vert
\begin{array}{l}
ik\\
j
\end{array}\right.,\left.\begin{array}{l}
ik\\
j
\end{array}\right>=q^{3}\left\vert
\begin{array}{l}
ik\\
j
\end{array}\right.,\left.\begin{array}{l}
ik\\
j
\end{array}\right>+(q^{5}-q^{-1})[2]^{-1/2}\left\vert
\begin{array}{l}
ik\\
k
\end{array}\right.,\left.\begin{array}{l}
ij\\
j
\end{array}\right> -(q^{4}-q^{-2})[2]^{-1/2}\left\vert
\begin{array}{l}
jk\\
k
\end{array}\right.,\left.\begin{array}{l}
ii\\
j
\end{array}\right>$$

$$R\left\vert
\begin{array}{l}
ik\\
k
\end{array}\right.,\left.\begin{array}{l}
ij\\
j
\end{array}\right>=q\left\vert
\begin{array}{l}
ij\\
j
\end{array}\right.,\left.\begin{array}{l}
ik\\
k
\end{array}\right>+(q^{3}-q){[3]^{1/2}\over{[2]}}\left(\left\vert
\begin{array}{l}
ij\\
k
\end{array}\right.,\left.\begin{array}{l}
ik\\
j
\end{array}\right> +\left\vert
\begin{array}{l}
ik\\
j
\end{array}\right.,\left.\begin{array}{l}
ij\\
k
\end{array}\right>\right)+$$

$$(q^{5}-q^{-1})[2]^{-1}\left\vert
\begin{array}{l}
ik\\
j
\end{array}\right.,\left.\begin{array}{l}
ik\\
j
\end{array}\right>+(q^{-1}-q-q^{3}+q^{5})\left\vert
\begin{array}{l}
ik\\
k
\end{array}\right.,\left.\begin{array}{l}
ij\\
j
\end{array}\right>$$

$$+(q^{-1}-q)\left\vert
\begin{array}{l}
jj\\
k
\end{array}\right.,\left.\begin{array}{l}
ii\\
k
\end{array}\right>-[2]^{2}\left\vert
\begin{array}{l}
jk\\
k
\end{array}\right.,\left.\begin{array}{l}
ii\\
j
\end{array}\right>$$

$$R\left\vert
\begin{array}{l}
jj\\
k
\end{array}\right.,\left.\begin{array}{l}
ii\\
k
\end{array}\right>=q\left\vert
\begin{array}{l}
ii\\
k
\end{array}\right.,\left.\begin{array}{l}
jj\\
k
\end{array}\right>+(q^{4}-1)\left\vert
\begin{array}{l}
ij\\
k
\end{array}\right.,\left.\begin{array}{l}
ij\\
k
\end{array}\right>+(q^{-1}-q)\left\vert
\begin{array}{l}
ik\\
k
\end{array}\right.,\left.\begin{array}{l}
ij\\
j
\end{array}\right>+$$

$$(q^{-1}-q-q^{3}+q^{5})\left\vert
\begin{array}{l}
jj\\
k
\end{array}\right.,\left.\begin{array}{l}
ii\\
k
\end{array}\right>-[2]^{2}\left\vert
\begin{array}{l}
jk\\
k
\end{array}\right.,\left.\begin{array}{l}
ii\\
j
\end{array}\right>,\eqno(42)$$
\vskip .3truecm
\noindent where we always assume that $1\leq i<j<k<l\leq n$. The above results
exhaust all the R matrix elements for $[1]\times [1]$, $[2]\times [2]$,
and $[1^{2}]\times [1^{2}]$. Because there are too many matrix elements
for $[21]\times [21]$, we only list those for $U_{q}(3)$ completely. But these
$U_{q}(3)$ R matrix elements are also those of $U_{q}(n)$. I. e., the results
are $n$ independent. Secondly, Another advantage of this method is that it
leaves coproduct multiplicity problem untouched in contrast to the projection
operator technique, in which one should derive CG matrix of the corresponding
quantum groups. Sometimes the coproduct concerned is not multiplicity-free. For
example, in our $[21]\times [21]$ case, the resultant irrep $[21]$
is of multiplicity two. However, in our method the R operator is acting
on the uncoupled basis vectors of the corresponding quantum groups, and has
nothing to do with the CG coefficients of $[21]\times [21]$.
From these simple examples, it can easily
be seen that the R matrix elements obtained in this way are $n$ independent.
I. e., all the R matrix elements of $U_{q}(n)$ listed apply for arbitrary $n$.
\vskip .3truecm
   The procedures for evaluating R matrix elements of other quantum groups
are very similar. For example, we can evaluate R matrix elements for the $G_{2}$
case with the help of Kalfagianni algebra$^{[22]}$. Hence the problem for constructing
R matrix is switched to find out all the possible algebraic
realizations of braid group. It should be noted that this method will also
become very tedious when the rank of the irrep increases. In such case,
the R operator is expressed in terms of a lengthy braid group generator product.
Though this method is not simpler than the other methods for higher dimensional
irreps of quantum group, it shows us a new way to calculate R matrix elements,
a new perspective of R matrices, and a transparent view of its braid group
structure.

\vskip 2truecm
\noindent {\bf Acknowledgement}
\vskip .4truecm
   The project is supported
by The Natural Science Foundation of China, and Excellent Young Teacher's
Foundation of The State Education Commission of China.
\vskip 1truecm
\noindent{\large\bf References}
\vskip .4truecm
\begin{tabbing}
\=1111\=2222222222222222222222222222222222222222222222222222222222222222222222\=\kill\\
\>{[1]}\>{Jimbo M, 1990 Yang-Baxter Equation in Integrable Systems, (Singapore, World Scientific). }\\
\>{[2]}\>{Yang C N, 1967 Phys. Rev. Lett. {\bf 19}, 1312}\\
\>{[3]}\>{Yang C N, Ge M L (eds.), 1989, Braid Group, Knot Theory and Statistical}\\
\>{}\>{Mechanics, (Singapore, World Scientific)}\\
\>{[4]}\>{T. Kohno, 1989  Braid Group, Knot Theory and Statistical Mechanics,}\\
\>{}\>{Yang C. N. and Ge M. L.  eds., (Singapore, World Scientific)}\\
\>{[5]}\>{Jones V F R, 1989  Australian National Univ. Rep. CMA-R23-89, 1989}\\
\>{[6]}\>{Ge M L et al, 1990 in Nonlinear Physics, Hu C H et al  eds. (Springer-Verlag, Berlin)}\\
\>{[7]}\>{Zhang R B, Gould M D, and Bracken A J, 1991 Nucl. Phys. {\bf B354}, 625}\\
\>{[8]}\>{Drinfeld V, 1985 Dokl. SSR {\bf 283}, 1060}\\
\>{[9]}\>{Jimbo M, 1986 Commun. Math. Phys. {\bf 102},537}\\
\>{[10]}\>{Jimbo M, 1985 Lett. Math. Phys. {\bf 10}, 63;{\bf 11}, 247}\\
\>{[11]}\>{Reshetikin N, 1987 LMOT Preprint E-4-87}\\
\>{[12]}\>{Akutsu Y, and Wadati M, 1986 J. Phys. Soc. Japan {\bf 55}, 1092;1466; 1880}\\
\>{[13]}\>{Schultz C L, 1981 Phys. Rev. Lett. {\bf 46},629}\\
\>{[14]}\>{Babelon O, de Vega H J and Viallet C M, 1981 Nucl. Phys. {\bf 190}, 542}\\
\>{[15]}\>{Witten E, 1989 Commun. Math. Phys. {\bf 121},351}\\
\>{[16]}\>{Ge M L, Li Y. Q and Xue K, 1990 Int. J. Mod. Phys. {\bf A5},1975}\\
\>{[17]}\>{Ge M L, Wu Y S and Xue K, 1991 Int. J. Mod. Phys. {\bf A6}, 3735}\\
\>{[18]}\>{Ge M L and Xue K, 1991 J. Math. Phys. {\bf 32}, 1991}\\
\>{[19]}\>{Yu Stroganov G, 1979 Phys. Lett. {\bf 74A}, 116}\\
\>{[20]}\>{Wenzl H, 1990 Commun. Math. Phys. {\bf 133}, 383}\\
\>{[21]}\>{Pan F. and Chen J Q, 1993 J. Math. Phys. {\bf 34},4305; 4316}\\
\>{[22]}\>{Kalfagianni E, 1993, J. Knot Theory and Its  Ramifications {\bf 2}, 431}\\
\end{tabbing}
\end{document}